# Pure L-Functions from Algebraic Geometry over Finite Fields

Daqing Wan

Department of Mathematics, University of California, Irvine, CA 92697-3875
dwan@math.uci.edu

**Abstract.** This survey gives a concrete and intuitively self-contained introduction to the theory of pure L-functions arising from a family of algebraic varieties defined over a finite field of characteristic $p$. The standard fundamental questions in any theory of L-functions include the meromorphic continuation, functional equation, Riemann hypothesis (RH for short), order of zeros at special points and their special values. Our emphasis here will be on the meromorphic continuation and the RH. These two questions can be described in a general setup without introducing highly technical terms.

The construction of the pure L-function depends on the choice of an absolute value of the rational number field $\mathbf{Q}$. In the case that the absolute value is the complex or $\ell$-adic absolute value ($\ell \neq p$), the rationality of the pure L-function requires the full strength of the $\ell$-adic cohomology including Deligne's main theorem on the Weil conjectures. In the case that the absolute value is the $p$-adic absolute value, the pure L-function is no longer rational but conjectured by Dwork to be $p$-adic meromorphic. This conjecture goes beyond all existing $p$-adic cohomology theories. Its truth opens up several new directions including a possible $p$-adic RH for such pure L-functions.

The guiding principle of our exposition in this paper is to describe all theorems and problems as simple as possible, directly in terms of zeta functions and L-functions without using cohomological terms. In the case that this is not easy to do so, we simply give an intuitive discussion and try to convey a little feeling. Along the way, a number of natural open questions and conjectures are raised, some of them may be accessible to certain extent but others may be somewhat wild due to the lack of sufficient evidences.

## 1 Introduction

The most basic question in number theory is to understand the integers. In particular, for a given integer $N$, we need to understand the absolute value $\|N\|$ for every absolute value $\|?\|$ on the rational number field $\mathbf{Q}$. For the complex absolute value, this is to determine

$$\|N\| = |N| = ?$$

For the $p$-adic absolute value with $p$ being a prime, this is to determine

$$|N|_p = p^{-a_p}, \ a_p = \mathrm{ord}_p(N) = ?$$

This last theoretical question is practically the problem of factoring integers which has important applications.

More generally, suppose that we are given a sequence of interesting integers
$$\{N_1, N_2, \cdots, \}$$
In order to understand this sequence of integers, one naturally forms a suitable generating function $Z(\{N_i\}, T)$ which contains all information about the given sequence. The basic question is then to understand the analytic properties of the generating function $Z(\{N_i\}, T)$ with respect to each absolute value $\|?\|$ of $\mathbf{Q}$. This includes the possible meromorphic continuation $Z(\{N_i\}, T)$ and a suitable RH about its zeros and poles, for both the complex absolute value and the $p$-adic absolute value. If we have a family of such generating functions, then we would like to understand its analytic variation when the parameter varies.

The most interesting type of sequences arises from counting prime ideals in a finitely generated commutative ring or equivalently from counting rational points on an algebraic variety. In the case of counting prime numbers in the ring $\mathbf{Z}$ of integers, the natural generating function is the Riemann zeta function. This first example was studied by Riemann from complex point of view and by Kummer-Kubota-Leopoldt from $p$-adic point of view. It is the motivating example for much of the modern developments on general Hasse-Weil zeta functions of algebraic varieties as well as their conjectural $p$-adic analogues.

Our interesting sequence of integers in this paper arises from counting rational points over various finite extension fields of an algebraic variety $X$ defined over a finite field of characteristic $p$. The resulting generating function is the zeta function of $X$ which is the object of study in the celebrated Weil conjectures. The zeta function is a rational function as proved by Dwork using $p$-adic methods. It satisfies a suitable complex and $\ell$-adic RH as proved by Deligne using $\ell$-adic methods, where $\ell$ is a prime number different from $p$. The $p$-adic RH for the zeta function is more complicated and remains mysterious in general. The variation of the whole zeta function, when the variety moves through an algebraic family, leads to new interesting questions which are understood to certain extent.

The zeta function is however not pure. That is, the zeros and poles have different absolute values. This is especially so from $p$-adic point of view. Thus, the zeta function decomposes as a product of pure pieces defined in terms of the absolute values of the zeros and poles. A finer form of the RH is to understand the purity decomposition. A further question is to understand the variation of each pure piece of the zeta function when the variety moves through an algebraic family. This naturally leads to the construction of pure L-functions arising from algebraic geometry. Our fundamental question here is then to understand the analytic properties of such a pure L-function, notably its meromorphic continuation and RH. Since the zeta function has integer coefficients, there are three different types of absolute values (complex, $\ell$-adic and $p$-adic) that we can choose to work. These lead to different results and different theories.

In the case that the absolute value is the complex or the $\ell$-adic absolute value, Deligne's main theorem shows that the pure L-function from algebraic geometry can be identified with a geometric L-function, that is the L-function of a certain geometric constructible $\ell$-adic étale sheaf. One can then apply the full machinery of $\ell$-adic étale cohomology. In particular, the pure L-function from algebraic geometry is always rational by Grothendieck's rationality theorem. It satisfies a suitable complex and $\ell$-adic RH by Deligne's theorem. The situation is quite algebraic in nature. All the expected finiteness properties hold. The main point is that one does not need to distinguish the subcategory of geometric $\ell$-adic sheaves from the full category of $\ell$-adic sheaves. Much of the relevant theory works for every $\ell$-adic sheaf, whether it is geometric or not.

In the case that the absolute value is the $p$-adic absolute value, the situation is quite different and much more complicated. A pure L-function from algebraic geometry is not rational any more. However, Dwork conjectured that a pure L-function from algebraic geometry is $p$-adic meromorphic. The situation is quite transcendental in nature. Grothendieck's specialization theorem, Katz's isogeny theorem and Berthelot's finiteness theorem on relative crystalline cohomology show that, at least in nice cases, a pure L-function from algebraic geometry can be identified with the L-function of a certain geometric $p$-adic étale sheaf. The trouble is that the L-function of a general $p$-adic étale sheaf does not behave well. The usual trace formula does not hold. Even worse, the L-function is not meromorphic in general, unlike what Katz conjectured. Thus, in order to prove Dwork's conjecture, one must distinguish the subcategory of geometric $p$-adic sheaves from the full category of $p$-adic sheaves. Our recent work shows that the geometric $p$-adic sheaves can be understood by introducing a new category with growth condition. Roughly speaking, this new category consists of infinite nuclear complexes of infinite rank nuclear overconvergent F-isocrystals. The nuclear overconvergent condition insures that the L-function is $p$-adic meromorphic. This then establishes the meromorphic continuation of a pure L-function from algebraic geometry and thus proves Dwork's conjecture. The $p$-adic RH for such a pure L-function is extremely mysterious. A good understanding seems to require entirely new ideas.

We would like to point out that there is a more general and more difficult type of zeta functions arising from counting algebraic cycles on an algebraic variety $X$ defined over a finite field. These zeta functions are called the zeta functions of algebraic cycles. They seem to be out of reach at this time. For zero cycles, they reduce to the zeta functions studied in the Weil conjectures, which are already quite interesting and fruitful. In general, they contain important arithmetic information about algebraic cycles and are related to Tate's conjecture. Under a mild finiteness condition on the effective cone of

the Chow group, these zeta functions of algebraic cycles are conjectured to be
$p$-adic meromorphic. If this conjecture is true, one could go on to understand
its $p$-adic RH, their variation when the variety moves through an algebraic
family, the purity decomposition and the resulting pure L-functions of algebraic cycles. Any proof of this meromorphic conjecture would likely have a
profound impact on arithmetic and geometry of algebraic cycles.

## 2  Rationality of zeta functions

Let $\mathbf{F}_q$ be the finite field of $q$ elements of characteristic $p$. Let $X$ be an algebraic variety defined over $\mathbf{F}_q$, namely, a separated scheme of finite type over
$\mathbf{F}_q$. For example, if $X$ is affine, then $X$ is defined by a system of polynomial
equations
$$f_1(x_1, \cdots, x_n) = \cdots = f_r(x_1, \cdots, x_n) = 0$$
in some affine $n$-space $\mathbf{A}^n$, where each $f_i$ is a polynomial defined over $\mathbf{F}_q$.
For an extension field $\mathbf{F}_{q^k}$ of degree $k$ over $\mathbf{F}_q$, let $X(\mathbf{F}_{q^k})$ denote the set of
$\mathbf{F}_{q^k}$-rational points on $X$. The zeta function of $X/\mathbf{F}_q$ is then defined to be
the following formal power series

$$Z(X/\mathbf{F}_q, T) = \exp(\sum_{k=1}^{\infty} \frac{\#X(\mathbf{F}_{q^k})}{k} T^k)$$

$$= \prod_{x \in X(\bar{\mathbf{F}}_q)} \frac{1}{(1 - T^{\deg(x)})^{1/\deg(x)}}$$

$$= \prod_{x \in X_0} \frac{1}{1 - T^{\deg(x)}} \in \mathbf{Z}[[T]],$$

where $X_0$ is the set of closed points on $X/\mathbf{F}_q$ and $\bar{\mathbf{F}}_q$ denotes a fixed algebraic
closure of $\mathbf{F}_q$. Recall that a closed point on $X/\mathbf{F}_q$ is simply the orbit of
an actual geometric point $x \in X(\bar{\mathbf{F}}_q)$ under the $q$-th power Frobenius map
$\sigma : x \to x^q$ and $\deg(x)$ is the smallest positive integer $k$ such that $\sigma^k(x) = x$.
The integer $\#X(\mathbf{F}_{q^k})$ is simply the number of fixed points of the $k$-th power
$\sigma^k$ acting on $X(\bar{\mathbf{F}}_q)$. We shall write $Z(X, T)$ for $Z(X/\mathbf{F}_q, T)$ when the ground
field $\mathbf{F}_q$ is clear.

The zeta function is a generating function for counting rational points on
the variety $X$ over various finite extension fields of $\mathbf{F}_q$. Although at each stage
$\mathbf{F}_{q^k}$, only finitely many points are counted, the generating function counted
all points of the variety $X$ over the algebraic closure $\bar{\mathbf{F}}_q$. This explains why
the zeta function should contain a great deal of geometric and arithmetic
information about the variety $X$. Thus, our first fundamental question is

*Question 2.1.* Understand the zeta function $Z(X, T)$.

A general principle in analytic arithmetic algebraic geometry is that all zeta functions and L-functions arising naturally from arithmetic algebraic geometry are analytically good functions. Here we are considering algebraic geometry over $\mathbf{F}_q$. Based on earlier results in various special cases such as diagonal hypersurfaces, curves and abelian varieties, Weil conjectured the following rationality result.

**Theorem 2.2.** *The zeta function $Z(X,T)$ is a rational function in $T$.*

This theorem was first proved by Dwork [Dw1] using $p$-adic analysis. His basic idea is to establish a trace formula which expresses the zeta function as a finite alternating product of the Fredholm determinants of several nuclear operators acting on certain $p$-adic Banach spaces. Since the Fredholm determinant of a nuclear operator is entire [Se], it follows that the zeta function is $p$-adic meromorphic. One then concludes the proof of Theorem 2.2 with the following $p$-adic analogue of Borel's rationality criteria.

**Lemma 2.3.** *Let $f(T) \in \mathbf{Z}[[T]]$ be a power series with integer coefficients. Then, $f(T)$ is a rational function if and only if $f(T)$ is analytic at the origin as a complex function and meromorphic as a $p$-adic function for some prime number $p$.*

Dwork's rationality proof pioneered his $p$-adic theory of zeta functions. Although his proof is not cohomological in nature, it can be viewed as a proof on the chain level in a suitable sense. Its refinement by introducing commuting differential operators is cohomological in nature. The relevant cohomology is of De Rham type as shown by Katz [K1]. This motivated later more systematic developments of various $p$-adic cohomology theories such as the formal cohomology [MW] for smooth affine varieties, the liftable cohomology [Lu] for smooth projective liftable varieties, crystalline cohomology for smooth projective varieties and rigid cohomology [Be] for arbitrary varieties. These $p$-adic cohomology theories in their full generality have not been well understood. The relevance of differential operators suggests the existence of a close relationship between zeta functions and differential equations. This is the subject of F-crystals [K3], which we will not discuss here.

Weil's conjectural approach to Theorem 2.2, generalizing his method for curves and abelian varieties, is to construct a suitable cohomology theory for varieties in characteristic $p$ so that the Lefschetz fixed point theorem holds which immediately implies the rationality of the zeta function. This led to the introduction and full development of the $\ell$-adic étale cohomology theory by Grothendieck and others, where $\ell$ is a prime number different from $p$. This theory provides some of the most powerful tools in the study of zeta functions, especially from the complex point of view. For the exceptional prime $p$, the theory of $p$-adic étale cohomology has also been studied by Grothendieck's school but it does not behave well for L-function purpose. Nevertheless, the $p$-adic étale cohomology contains important $p$-adic information about the zeta function. This turns out to be closely related to Dwork's conjecture. We shall however not discuss this cohomological point of view here due to our restricted simple nature of this survey.

The zeta function is a power series with integer coefficients. However, all proofs of Theorem 2.2 use non-archimedian methods, either $\ell$-adic or $p$-adic. No direct "motivic" proof over $\mathbf{Q}$ or $\mathbf{C}$ is known. Although the classical method of Gauss sums and Jacobi sums for the elementary diagonal case can be viewed as a method over $\mathbf{C}$, it seems hopeless to extend this method to the general case. In the case of curves, there is a direct approach over $\mathbf{Z}$ using the Riemann-Roch theorem. It would be of great interest to extend this Riemann-Roch approach to the general case. Such a direct approach would make it possible to attack Tate's conjecture [Ta] relating orders of zeros and ranks of algebraic cycles which so far has been inaccessible. Such a Riemann-Roch approach should also be useful in understanding the harder zeta functions of algebraic cycles introduced in [W1].

The first application of Theorem 2.2 is the existence of a formula for the number of rational points on $X$. Since the zeta function is rational, there are finitely many algebraic integers $\alpha_i$ and $\beta_j$ such that

$$Z(X, T) = \frac{\prod_i (1 - \alpha_i T)}{\prod_j (1 - \beta_j T)}.$$

Taking logarithmic derivative, one deduces the following well structured formula for every positive integer $k$:

$$\#X(\mathbf{F}_{q^k}) = \beta_1^k + \beta_2^k + \cdots - \alpha_1^k - \alpha_2^k - \cdots.$$

In particular, this provides a fast algorithm to compute the number of rational points on the variety $X$ over a large finite field $\mathbf{F}_{q^k}$ provided $X$ is defined over a small finite field $\mathbf{F}_q$. The theory can be improved to get a fast algorithm, even for $X$ defined over a large finite field $\mathbf{F}_q$ as long as the characteristic $p$ is small ($q$ can be large), see [W5] for a perspective on this algorithmic subject which has important practical applications.

Once we know that the zeta function is a rational function, we can move on to the next fundamental question. What can we say about its zeros and poles? Ideally, we would like to know how many zeros and poles with a given absolute value. This is the RH for the zeta function. Since the reciprocal zeros $\alpha_i$ and $\beta_j$ are algebraic integers, there are several different types of absolute values (complex, $\ell$-adic and $p$-adic) that we can consider. Accordingly, we can talk about the complex RH, the $\ell$-adic RH and the $p$-adic RH for the zeta function $Z(X, T)$. These questions and their family versions are discussed in the following sections.

## 3 Purity decomposition and RH

The zeta function has rational coefficients. In order to understand its analytic properties, we have to choose an absolute value of the rational number field

**Q**. For this purpose, we let $\|?\|$ be a fixed absolute value on **Q**. Let $\Omega$ be the smallest extension field of **Q** such that $\Omega$ is both algebraically closed and topologically complete with respect to $\|?\|$. If $\|?\|$ is the complex absolute value $|?|$, $\Omega$ is the field **C** of complex numbers. If $\|?\|$ is the $p$-adic absolute value $|?|_p$ for some prime number $p$, $\Omega$ is the field $\mathbf{C}_p$ of $p$-adic numbers which is the completion of an algebraic closure of the $p$-adic rational numbers $\mathbf{Q}_p$, where the $p$-adic absolute value is normalized by $|p|_p = 1/p$. For a fixed prime power $q$ of $p$, we define the slope function on $\Omega$ by

$$s(\alpha) = \log_q \|\alpha\|, \text{ if } \|?\| = |?|$$

and

$$s(\alpha) = -\log_q \|\alpha\|, \text{ if } \|?\| \neq |?|.$$

This definition depends on $q$ which is our base of the slope function. For $\alpha \in \Omega$, we can write

$$\alpha = q^{s(\alpha)} u(\alpha),$$

where $u(\alpha)$ is a number with absolute value 1. Note that in the complex absolute value case, our slope function $s(\alpha)$ is twice the weight function $w(\alpha)$ as defined by Deligne. In the $p$-adic case, the slope function is simply the order function

$$s(\alpha) = \mathrm{ord}_q(\alpha).$$

For a given polynomial

$$P(T) = \prod_i (1 - \alpha_i T) \in \Omega[T],$$

and a real number $s$, we define the slope $s$ part of $P(T)$ by

$$P_s(T) = \prod_{s(\alpha_i)=s} (1 - \alpha_i T) \in \Omega[T].$$

This immediately yields the purity (or slope) decomposition of $P(T)$:

$$P(T) = \prod_{s \in \mathbf{R}} P_s(T).$$

This is a finite product since $P(T)$ is a polynomial. This definition easily extends to rational functions in $\Omega(T)$ as well as meromorphic functions in $\Omega((T))$ by the Weierstrass factorization theorem. In the latter case, the purity decomposition is an infinite product in general. Applying the purity decomposition to the zeta function $Z(X,T)$, our notation becomes

$$Z(X,T) = \prod_{s \in \mathbf{R}} Z_s(X,T),$$

where $Z_s(X,T)$ is the slope $s$ part of $Z(X,T)$. This is called the purity (or slope) decomposition of the zeta function $Z(X,T)$. Our next fundamental question is then to understand this purity decomposition. That is,

*Question 3.1.* Understand each pure slope $s$ part $Z_s(X,T)$ of $Z(X,T)$.

The first step is to understand the degree of the rational function $Z_s(X,T)$ for each $s$. To be precise, we recall that the degree of a rational function is the degree of the numerator minus the degree of the denominator. Similarly, the total degree of a rational function is the degree of the numerator plus the degree of the denominator.

**Definition 3.2.** Let $d(X)$ (resp. $D(X)$) denote the degree (resp. the total degree) of the zeta function $Z(X,T)$. Similarly, for each real number $s$, let $d_s(X)$ (resp. $D_s(X)$) denote the degree (resp. the total degree) of the slope $s$ part $Z_s(X,T)$ of the zeta function $Z(X,T)$.

It is clear that we have the purity decomposition for the degrees $d(X)$ and $D(X)$:

$$d(X) = \sum_{s \in \mathbf{R}} d_s(X),$$

$$D(X) = \sum_{s \in \mathbf{R}} D_s(X).$$

The RH for the zeta function is to determine the exact slopes of the zeros and poles. It is easy to see that the number of reciprocal zeros of slope $s$ is given by $(D_s(X) + d_s(X))/2$. Similarly, the number of reciprocal poles of slope $s$ is given by $(D_s(X) - d_s(X))/2$. Thus, the following weaker but more precise form of Question 3.1 is already the RH for $Z(X,T)$.

*Question 3.3.* Understand the pure degree $d_s(X)$ and the pure total degree $D_s(X)$ for all $s$.

Both the degree $d(X)$ and the total degree $D(X)$ of the whole zeta function $Z(X,T)$ can be effectively bounded using $p$-adic methods as shown by Bombieri [Bo]. This is because Dwork's $p$-adic theory is constructive. Thus, the pure degree $d_s(X)$ and the pure total degree $D_s(X)$ are also effectively bounded for all $s$. The integers $d_s(X)$ and $D_s(X)$ depend both on the slope $s$ and on the variety $X$. Of course, they also depend on the choice of the absolute value $\|?\|$ which was built in the definition of the slope decomposition. In this section, we consider the case that $X$ is fixed. In next section, we consider how $d_s(X)$ and $D_s(X)$ vary when $X$ varies.

Deligne's main theorem [De2] on the complex and the $\ell$-adic RH can be stated in our notations as follows.

**Theorem 3.4 (complex case).** *Let $\|?\|$ be the complex absolute value. Let $n$ be the dimension of the variety $X$. If*

$$s \notin \{0, \frac{1}{2}, 1, \frac{3}{2}, 2, \cdots, n\},$$

*then*

$$d_s(X) = D_s(X) = 0, \ Z_s(X,T) = 1.$$

This shows that for the complex absolute value, the non-trivial slopes are rational numbers in the interval $[0, n]$ with denominators at most 2. Thus, for the complex RH, it remains to determine the $2n + 1$ values $d_s(X)$ and $D_s(X)$, where $s$ varies in the above exceptional set of $2n+1$ numbers. These remaining values are in general difficult to determine although some extremal cases such as $s = n, n-\frac{1}{2}$ can be done. They depend on the detailed geometry of the variety $X$. However, in nice situations, they can be determined by the Betti numbers. This includes the smooth projective (more generally smooth proper) case as conjectured by Weil and first proved by Deligne [De1] [De2] using $\ell$-adic cohomology. A similar proof was later given by Faltings [Fa] using crystalline cohomology.

For the $\ell$-adic RH, the answer is much simpler and very clean.

**Theorem 3.5 ($\ell$-adic case).** *Let $\|?\|$ be the $\ell$-adic absolute value for some prime $\ell \neq p$. If $s \neq 0$, then*

$$d_s(X) = D_s(X) = 0, \; Z_s(X, T) = 1.$$

*In particular,*

$$D_0(X) = D(X), \; d_0(X) = d(X), \; Z_0(X, T) = Z(X, T).$$

*That is, all zeros and poles of the zeta function are $\ell$-adic units.*

Since the $\ell$-adic case is always pure and thus gives no interesting decomposition, from now on, we shall mostly restrict our attention to the complex case and the $p$-adic case.

For the $p$-adic RH, unfortunately, no clean general answer is possible, even in the smooth projective case, even in the case of smooth projective curves. However, one has the following weak but simple general result, which is a consequence of the rationality of the zeta function.

**Theorem 3.6 ($p$-adic case).** *Let $\|?\|$ be the $p$-adic absolute value. There is an effectively computable positive integer $N$ such that if*

$$s \notin \{0, \frac{1}{N}, \frac{2}{N}, \cdots, n\},$$

*then*

$$d_s(X) = D_s(X) = 0, \; Z_s(X, T) = 1,$$

*where $n$ is the dimension of $X$.*

This result is far weaker than Theorem 3.4, because the denominator $N$ here is not bounded by 2. In fact, the denominator $N$ cannot be bounded by any finite absolute constant. It depends very much on the variety $X$ and the prime number $p$, not just on the geometry of $X$. This explains why the $p$-adic RH for the zeta function is very complicated. It can be determined in a few special cases such as the elementary diagonal hypersurface case where one can use the Stickelberger theorem for Gauss sums. In nice situations such as the smooth projective case, a good lower bound for the Newton polygon (which determines the $p$-adic RH) is given by the Hodge polygon (constructed using the Hodge numbers of a lifting of $X$), as conjectured by Katz and proved by Mazur [M2]. Strictly speaking, both polygons are defined for each cohomological dimension. For our restricted purpose, one could either fix a cohomological dimension or take the collection over all cohomological dimensions. If the two polygons coincide, the variety $X$ is called ordinary.

For a given smooth projective variety $X$, Mazur's theorem provides a geometric lower bound for the arithmetic Newton polygon but it does not tell if $X$ is ordinary or how far $X$ is from being ordinary. There is no known clean recipe (conjectural or not) to determine the Newton polygon of $X$. A preliminary step might be to look at the size of the endomorphism group of $X$. The larger the endomorphism group for $X$ is, the more relations there would be among the zeros and poles and thus $X$ would be less likely ordinary, as one observes in the diagonal case and the supersingular elliptic curve case. It would be interesting to make such heuristic arguments more precise.

Another interesting question is to consider the limiting behavior of the Newton polygon as $p$ varies. For instance, let $X$ be a smooth projective variety defined over $\mathbf{Q}$ and let $X_p$ be the reduction of $X$ mod $p$ for large prime $p$. Let $NP(X_p)$ denote the Newton polygon of $X_p$ and let $HP(X)$ denote the Hodge polygon of $X$. As $p$ goes to infinity, the Newton polygon $NP(X_p)$ would not have a limit in general. But it always has its lower limit. As this lower limit has nothing to do with any particular prime $p$, it should be a geometric invariant of $X$. Thus, it is tempting to make

**Conjecture 3.7.** *Let $X$ be a smooth projective variety over $\mathbf{Q}$. Then,*

$$\lim_{p\to\infty} \inf NP(X_p) = HP(X).$$

Since there are only finitely many possibilities for the $NP(X_p)$ for a fixed smooth projective $X/\mathbf{Q}$, the above conjecture is equivalent to saying that there are infinitely many ordinary primes $p$ (probably of positive density) for a given smooth projective $X/\mathbf{Q}$. If the endomorphism group of $X$ is sufficiently small, one may further hope that the set of ordinary primes $p$ for $X$ often has density 1. A somewhat related conjecture is given by Serre [Oo] in the case of abelian varieties over number fields. Similarly, one could ask

$$\lim_{p\to\infty} \sup NP(X_p) = ?$$

This sup limit should again be a geometric invariant of $X$.

So far, we have been concerned with the degree of the pure slope $s$ part $Z_s(X,T)$. A finer question is to understand the rational function $Z_s(X,T)$ itself for non-trivial slope $s$. For instance, one could ask about the possible rationality of the coefficients for $Z_s(X,T)$ when $Z_s(X,T)$ is written as the quotient of two relatively prime polynomials with constant term 1. The answer depends on the absolute value $\|?\|$ we choose. If $\|?\|$ is the $\ell$-adic absolute value, then Theorem 3.5 shows that $Z_s(X,T)$ trivially has integer coefficients. The same result holds in the complex absolute value case. This follows from Deligne's main theorem and Galois theory.

**Theorem 3.8.** *Let $\|?\|$ be the complex absolute value. Then, the pure slope $s$ part $Z_s(X,T)$ has integer coefficients.*

In the $p$-adic case, it is no longer true that the coefficients of $Z_s(X,T)$ are integers, because the purity decomposition becomes more substantial, even in nice situations. One has the following $p$-adic rationality result which is a consequence of the rationality of $Z(X,T)$.

**Theorem 3.9.** *Let $\|?\|$ be the p-adic absolute value. Then, the coefficients of the pure slope $s$ part $Z_s(X,T)$ are p-adic integers in $\mathbf{Z}_p$, which are also algebraic integers.*

Of course, the reciprocal roots of $Z_s(X,T)$ will not be in $\mathbf{Z}_p$ in general.

## 4 Variation of the pure degrees

In the previous section, we discussed the pure degrees $d_s(X)$ and $D_s(X)$ for the zeta function of a single variety $X$. In this section, we turn to discussing how the pure degrees $d_s(X)$ and $D_s(X)$ vary when $X$ moves through an algebraic family.

Let $f: Y \to X$ be a family of algebraic varieties over $\mathbf{F}_q$ parametrized by $X$. For each geometric point $x \in X(\mathbf{F}_{q^{\deg(x)}})$, the fibre $Y_x = f^{-1}(x)$ is an algebraic variety defined over $\mathbf{F}_{q^{\deg(x)}}$ and thus we have the purity decomposition:

$$Z(Y_x/\mathbf{F}_{q^{\deg(x)}}, T) = \prod_{s \in \mathbf{R}} Z_s(Y_x/\mathbf{F}_{q^{\deg(x)}}, T),$$

where the slope function is defined with respect to $q^{\deg(x)}$ since $\mathbf{F}_{q^{\deg(x)}}$ is the base field of $Y_x$. For each rational number $s$, let $d_s(Y_x)$ (resp. $D_s(Y_x)$) denote the degree (resp. the total degree) of the slope $s$ part $Z_s(Y_x, T)$ of the zeta function $Z(Y_x, T)$, where we stress again that the slope function is defined with respect to $q^{\deg(x)}$. The general results in Section 3 give information about the nature of the numbers $d_s(Y_x)$ and $D_s(Y_x)$ for a fixed $x$. We would like to understand how these numbers $d_s(Y_x)$ and $D_s(Y_x)$ vary when the geometric point $x$ varies. For a given integer $m$, we define

$$X(d_s, m) = \{x \in X | d_s(Y_x) = m\}, \ X(D_s, m) = \{x \in X | D_s(Y_x) = m\}.$$

These are subsets of $X$. We would like to understand the possible algebraic and geometric structure of these sets.

The theory of $\ell$-adic cohomology and Deligne's main theorem [De2] imply

**Theorem 4.1.** *Let $f : Y \to X$ be a family of algebraic varieties over $\mathbf{F}_q$. Let $\|?\|$ be the complex or the $\ell$-adic absolute value. Then for every integer $m$ and every slope $s$, the set $X(d_s, m)$ is a constructible subset of $X$.*

The same general statement for the harder $X(D_s, m)$ does not seem to follow from the existing results because the nature of cancellation is not well understood in general. Thus, we have

*Question 4.2.* Let $f : Y \to X$ be a family of algebraic varieties over $\mathbf{F}_q$. Let $\|?\|$ be the complex or the $\ell$-adic absolute value. Is it true that for every integer $m$ and every slope $s$, the set $X(D_s, m)$ is a constructible subset of $X$?

This is known to be true if $f$ is smooth and proper because there is then no cancellation of zeros and poles by Deligne's theorem. We hope and are slightly inclined to believe that Question 4.2 has a positive answer, at least after replacing the constructible notion by the more general definable notion of logic.

In the $p$-adic case, one can ask a similar question.

**Conjecture 4.3.** *Let $f : Y \to X$ be a family of algebraic varieties over $\mathbf{F}_q$. Let $\|?\|$ be the $p$-adic absolute value. Then for every integer $m$ and every slope $s$, the set $X(d_s, m)$ is a constructible subset of $X$.*

*Question 4.4.* Let $f : Y \to X$ be a family of algebraic varieties over $\mathbf{F}_q$. Let $\|?\|$ be the $p$-adic absolute value. Is it true that for every integer $m$ and every slope $s$, the set $X(D_s, m)$ is a constructible (or more generally definable) subset of $X$?

Conjecture 4.3 and Question 4.4 are known to have a positive answer in nice situations such as when $f : Y \to X$ is proper and smooth over $\mathbf{F}_q$ with a proper smooth lifting to characteristic zero. In such a case, Conjecture 4.3 is a consequence of Grothendieck's specialization theorem [K4] and Berthelot's finiteness theorem [Be] on the relative crystalline cohomology or the relative rigid cohomology in the proper smooth liftable case. The stronger Question 4.4 in such a nice situation also needs Deligne's theorem to avoid the cancellation problem.

The general case of Conjecture 4.3 follows from the conjectural finiteness of the relative rigid cohomology. However, Conjecture 4.3 is already known to be true for $s = 0$ by the congruence formula of Deligne-Katz [De3] [K2] for zeta functions. The general case of Question 4.4 is apparently more difficult. It is not known to be true even for $s = 0$.

In some nice cases such as the universal family of hypersurfaces (or more generally complete intersections), the generic Newton polygon coincides with the Hodge polygon as conjectured by Mazur [M1]. This can be proved in two ways. One approach is to use hyperplane sections to reduce the question to the case of a generic plane curve, as worked out by Illusie [Il] using crystalline cohomology and some ideas of Deligne. Another more flexible approach introduced by the author [W2] is to establish suitable local to global decomposition theorems to reduce the question to the diagonal case where Stickelberger theorem applies, see [W7] for an exposition of this method. This method has other applications such as the more general Adolphson-Sperber conjecture [AS] for the generic Newton polygon of exponential sums. Hence, in such a generic ordinary case, the generic values of $d_s(X)$ and $D_s(X)$ are determined by the Hodge numbers of $X$. Of course, for a given smooth projective hypersurface, there is still no simple recipe to determine if $X$ is ordinary.

It seems very difficult and complicated to have a complete understanding of the stratification of the universal family of hypersurfaces by Newton polygons. This is so, even in the case of curves. However, for the more managable family of abelian varieties, the stratification question by Newton polygons is reasonably well understood by the work of de Jong and Oort [DO].

## 5 Variation of the zeta function

Let $f : Y \to X$ be a family of algebraic varieties over $\mathbf{F}_q$ parametrized by $X$. In this section, we consider how the zeta function $Z(Y_x, T)$ varies when the parameter $x$ varies. A standard procedure is to understand all the higher moments of the zeros and poles of the family of rational functions $Z(Y_x, T)$. Write

$$Z(Y_x, T) = \frac{\prod_i (1 - \alpha_i(x) T)}{\prod_j (1 - \beta_j(x) T)},$$

where the $\alpha_i(x)$ and the $\beta_j(x)$ are algebraic integers parametrized by $x$. For a positive integer $k$, the $k$-th moment of the zeros and poles of $Z(Y_x, T)$ for $x \in X(\mathbf{F}_q)$ is given by the sum

$$S_k(f) = \sum_{x \in X(\mathbf{F}_q)} (\alpha_1(x)^k + \alpha_2(x)^k + \cdots - \beta_1(x)^k - \beta_2(x)^k - \cdots).$$

More generally, the $k$-th moment over the $d$-th extension field $\mathbf{F}_{q^d}$ of $\mathbf{F}_q$ is defined by

$$S_{k,d}(f) = \sum_{x \in X(\mathbf{F}_{q^d})} (\alpha_1(x)^{k \frac{d}{\deg(x)}} + \alpha_2(x)^{k \frac{d}{\deg(x)}} + \cdots$$
$$- \beta_1(x)^{k \frac{d}{\deg(x)}} - \beta_2(x)^{k \frac{d}{\deg(x)}} - \cdots).$$

This is an integer by Galois theory. The variation question is then to understand the $k$-th moment sequence $S_{k,d}(f)$ ($d = 1, 2, \cdots$) for every $k$. In terms of generating functions, we need to understand the $k$-th power L-function of the family $f$ defined by

$$L^{[k]}(f, T) = \exp(\sum_{d=1}^{\infty} \frac{S_{k,d}(f)}{d} T^d)$$

$$= \prod_{x \in X(\bar{\mathbf{F}}_q)} \frac{1}{Z(Y_x/\mathbf{F}_{q^{k\deg(x)}}, T^{\deg(x)})^{1/\deg(x)}}$$

$$= \prod_{x \in X_0} \frac{1}{Z(Y_x/\mathbf{F}_{q^{k\deg(x)}}, T^{\deg(x)})} \in \mathbf{Z}[[T]].$$

In the last two equations, the degree of $x$ is defined over $\mathbf{F}_q$, but the ground field for the fibre $Y_x$ has been extended from $\mathbf{F}_{q^{\deg(x)}}$ to its $k$-th extension field $\mathbf{F}^{k\deg(x)}$. Explicitly,

$$L^{[k]}(f, T) = \prod_{x \in X_0} \frac{\prod_j (1 - \beta_j^k(x) T^{\deg(x)})}{\prod_i (1 - \alpha_i^k(x) T^{\deg(x)})} \in \mathbf{Z}[[T]].$$

**Theorem 5.1.** *The $k$-th power L-function $L^{[k]}(f, T)$ of the family $f$ is a rational function for every positive integer $k$.*

This result is a special case of a more general rationality result for certain partial zeta functions studied in [W11]. It can be proved using either Dwork's $p$-adic method or Grothendieck's $\ell$-adic method. The only new ingredient is to use tensor operations and Newton's formula expressing the $k$-th power symmetric functions in terms of elementary symmetric functions. A theorem of Faltings [W11] shows that the partial zeta function is always nearly rational.

By Theorem 5.1, for each $k$, there are finitely many algebraic integers $\gamma_i(k)$ and $\delta_j(k)$ such that

$$L^{[k]}(f, T) = \frac{\prod_i (1 - \gamma_i(k) T)}{\prod_j (1 - \delta_j(k) T)}.$$

Thus, for every positive integer $k$, we have the formula for the $k$-th moment sequence

$$S_{k,d}(f) = \delta_1(k)^d + \delta_2(k)^d + \cdots - \gamma_1(k)^d - \gamma_2(k)^d - \cdots, \quad d = 1, 2, \cdots.$$

Such formulas imply the existence of a suitable equi-distribution theorem for the zeros and poles of the family of rational functions $Z(Y_x, T)$ parametrized by $x \in X$. The simple approach we take here is related to but different from the approach of Deligne-Katz [K5] via monodromy representations.

To be more precise in the application of the equi-distribution problem, one needs to know information about the RH for the $k$-th power L-function $L^{[k]}(f, T)$. Applying the purity decomposition to $L^{[k]}(f, T)$, we can write

$$L^{[k]}(f, T) = \prod_{s \in \mathbf{R}} L_s^{[k]}(f, T),$$

where $L_s^{[k]}(f, T)$ is the slope $s$ part of $L^{[k]}(f, T)$ and the slope is defined with respect to the base $q^k$. Recall that the purity decomposition always depends on our choice of the absolute value $\|?\|$ on $\mathbf{Q}$.

**Definition 5.2.** Let $d(f,k)$ (resp. $d_s(f,k)$) be the degree of the rational function $L^{[k]}(f,T)$ (resp. $L_s^{[k]}(f,T)$). Similarly, let $D(f,k)$ (resp. $D_s(f,k)$) be the total degree of the rational function $L^{[k]}(f,T)$ (resp. $L_s^{[k]}(f,T)$).

Although the present situation is more complicated than the zeta function case treated in Section 3, it is essentially of the same nature. That is, similar RH holds for $L^{[k]}(f,T)$ in the complex and $\ell$-adic case. One simply keeps track of the simple proof of Theorem 5.1 and applies Deligne's theorem. In particular, the zeros and poles of $L^{[k]}(f,T)$ are always $\ell$-adic units. For the complex absolute value, we make it explicit as follows, see [W11] for a more general result on partial zeta functions.

**Theorem 5.3 (complex case).** *Let $\|?\|$ be the complex absolute value. Let $n$ be the dimension of the variety $X$ and let $m$ be the relative dimension of the family $f$. For each non-negative integer $k$, if*

$$s \notin \{0, \frac{1}{2}, 1, \frac{3}{2}, 2, \cdots, km+n\},$$

*then*

$$d_s(f,k) = D_s(f,k) = 0, \ L_s^{[k]}(f,T) = 1.$$

This result shows that it suffices to understand $d_s(f,k)$ and $D_s(f,k)$ for the above $2(km+n)+1$ exceptional values of the slope $s$. These are in general quite complicated to determine except in some extremal cases.

As for the $p$-adic RH, one only has the following very weak general result.

**Theorem 5.4 ($p$-adic case).** *Let $\|?\|$ be the $p$-adic absolute value. There is an effectively computable positive integer $N(f,k)$ such that if*

$$s \notin \{0, \frac{1}{N(f,k)}, \frac{2}{N(f,k)}, \cdots, km+n\},$$

*then*

$$d_s(f,k) = D_s(f,k) = 0, \ L_s^{[k]}(f,T) = 1,$$

*where $n$ is the dimension of $X$ and $m$ is the relative dimension of $f$.*

Again, this result for the $p$-adic case is far weaker than Theorem 5.3 for the complex case. The problem is that we know almost nothing about the denominator $N(f,k)$. An important problem is to estimate the size of $N(f,k)$. This is related to the size of the degree $d_s(f,k)$. The extra integer $k$ provides a new dimension of variation. That is, how the degrees $d(f,k)$, $D(f,k)$, $d_s(f,k)$ and $D_s(f,k)$ vary with $k$. Although the size of the degrees $d(f,k)$ and $D_s(f,k)$ can be effectively bounded, no explicit general bounds are available in the literature. This should be a realistic interesting problem to study.

The first example of great arithmetic interest is the universal family $f$ of elliptic curves parametrized by a modular curve or an Igusa curve. In this case, the sequence of integers $d(f,k)$ is mutually determined by the sequence of dimensions of modular forms of weight $k+2$. The harder sequence $D(f,k)$ is unknown because of possible cancellation of zeros and poles. For the complex absolute value, the easier degree sequence $d_s(f,k)$ is understood by Deligne's theorem on the Ramanujan-Peterson conjecture [De4] but again the harder total degree sequence $D_s(f,k)$ is unknown. For the $p$-adic absolute value, the easier degree sequence $d_s(f,k)$ is already quite mysterious and one conjectures [W6] that $N(f,k)$ is bounded independent of $k$, or even stronger, $d_s(f,k)$ is bounded independent of $k$ and $s$. This is the $p$-adic analogue of the Ramanujan-Peterson conjecture. Very little is known about it. The variation of $d_s(f,k)$ with $k$ is closely related to the Gouvêa-Mazur conjecture [GM], see Coleman [Co] and [W4] for positive results in this direction. The variation of $d_s(f,k)$ with $k$ is also closely related to the geometry of the eigencurve as studied by Coleman-Mazur [CM].

## 6 Variation of the pure part of the zeta function

Let $f: Y \to X$ be a family of algebraic varieties over $\mathbf{F}_q$ parametrized by $X$. In this section, for each fixed slope $s$, we consider how the pure part $Z_s(Y_x, T)$ varies when the parameter $x$ varies. This is the full aspect of Question 3.1. This question is deeper than just how the degree $d_s(Y_x)$ of $Z_s(Y_x, T)$ varies with $x$. It is also deeper than how the total zeta function $Z(Y_x, T)$ varies with $x$. As in Section 5, the standard procedure is to understand the $k$-th moment sequence associated to the reciprocal zeros and reciprocal poles of $Z_s(Y_x, T)$. Equivalently, we need to understand the $k$-th power L-function associated to the $k$-th moment sequence.

Write
$$Z_s(Y_x/\mathbf{F}_{q^{\deg(x)}}, T) = \frac{\prod_i (1 - \alpha_i(x,s)T)}{\prod_j (1 - \beta_j(x,s)T)}.$$

The associated $k$-th moment over the $d$-th extension field $\mathbf{F}_{q^d}$ of $\mathbf{F}_q$ is defined by

$$S_{k,d}(s,f) = \sum_{x \in X(\mathbf{F}_{q^d})} (\alpha_1(x,s)^{k\frac{d}{\deg(x)}} + \alpha_2(x,s)^{k\frac{d}{\deg(x)}} + \cdots$$
$$- \beta_1(x,s)^{k\frac{d}{\deg(x,s)}} - \beta_2(x,s)^{k\frac{d}{\deg(x)}} - \cdots).$$

This is an integer in the complex absolute value case by a generalization of Theorem 3.8. It is a $p$-adic integer in $\mathbf{Z}_p$ in the $p$-adic case by a generalization of Theorem 3.9. The variation question is then to understand the $k$-th moment slope $s$ sequence $S_{k,d}(s,f)$ $(d = 1, 2, \cdots)$ for every $k$. In terms of generating functions, we need to understand the following pure L-function from algebraic geometry.

**Definition 6.1.** Let $k$ be a positive integer. For a given slope $s \in \mathbf{Q}$, the $k$-th power slope $s$ L-function $L^{[k]}(s, f, T)$ attached to the family $f$ is defined to be

$$L^{[k]}(s, f, T) = \exp(\sum_{d=1}^{\infty} \frac{S_{k,d}(s, f)}{d} T^d)$$

$$= \prod_{x \in X(\bar{\mathbf{F}}_q)} \frac{1}{Z_s(Y_x/\mathbf{F}_{q^{k\deg(x)}}, T^{\deg(x)})^{1/\deg(x)}}$$

$$= \prod_{x \in X_0} \frac{1}{Z_s(Y_x/\mathbf{F}_{q^{k\deg(x)}}, T^{\deg(x)})} \in \Omega[[T]],$$

where the slope function is defined with respect to the base $q^k$ (note that we have already replaced $T$ by $T^{\deg(x)}$ in the Euler factors).

In the last two equations, the degree of $x$ is defined over $\mathbf{F}_q$, but the ground field for the fibre $Y_x$ has been extended from $\mathbf{F}_{q^{\deg(x)}}$ to its $k$-th extension field $\mathbf{F}_q^{k\deg(x)}$. Explicitly,

$$L^{[k]}(s, f, T) = \prod_{x \in X_0} \frac{\prod_j (1 - \beta_j^k(x, s) T^{\deg(x)})}{\prod_i (1 - \alpha_i^k(x, s) T^{\deg(x)})} \in \Omega[[T]].$$

The $k$-the power slope $s$ L-function arises in a natural way from arithmetic and geometry. Thus, by the general principle, we expect it to be an analytically good function. For the $\ell$-adic absolute value, the situation reduces to the $k$-th power L-function of Section 5, since the purity decomposition is trivial. For the complex absolute value, the situation is deeper than Theorem 5.1. One needs the full strength of Deligne's main theorem [De2] which says that the higher direct image sheaf is mixed. This together with Grothendieck's rationality theorem and Newton's formula implies the following result which can be viewed as the archimedian analogue of Dwork's conjecture.

**Theorem 6.2 (complex case).** *Let $\|?\|$ be the complex absolute value. For every positive integer $k$ and every rational number $s$, the $k$-th power slope $s$ L-function $L^{[k]}(s, f, T)$ is a rational function.*

In the case that the absolute value is the $p$-adic absolute value, the situation is much more complicated. For slope $s = 0$, the pure L-function $L^{[k]}(0, f, T)$ can be obtained as a suitable limit of the easier mixed L-functions treated in the previous section. To see this, one observes that the total degree of the zeta function of each fibre $Y_x$ is uniformly bounded by the results in [Bo]. This implies that for certain choices of positive integer $M$ depending on $f$, we have the following congruence

$$L^{[k+p^m M]}(f, T) \equiv L^{[k+p^{m+1} M]}(f, T) \pmod{p^m}$$

for all positive integers $m$. This $p$-adic continuity relation shows that the $p$-adic limit
$$\lim_{m\to\infty} L^{[k+p^m M]}(f,T)$$
exists as a formal $p$-adic power series. In fact, one checks easily from the definition of the pure L-function that the above $p$-adic limit is precisely the pure slope 0 L-function:
$$L^{[k]}(0,f,T) = \lim_{m\to\infty} L^{[k+p^m M]}(f,T).$$

For slope $s > 0$, we do not know any $p$-adic limiting formula for the pure slope L-function $L^{[k]}(s,f,T)$ in terms of the mixed L-functions $L^{[k]}(f,T)$. For slope $s = 0$, although each mixed L-function $L^{[k]}(f,T)$ is rational by Theorem 5.1, we cannot conclude that its limit $L^{[k]}(0,f,T)$ would be rational or even meromorphic. There are geometric examples (the universal family of elliptic curves, for instance) which show that the $k$-th power slope 0 L-function $L^{[k]}(0,f,T)$ is not rational in the $p$-adic case. Dwork, however, conjectured [Dw5] that the pure slope L-function is always $p$-adic meromorphic.

**Conjecture 6.3 (Dwork).** *Let $\|?\|$ be the $p$-adic absolute value. For every positive integer $k$ and every rational number $s$, the $k$-th power slope $s$ L-function $L^{[k]}(s,f,T)$ is a $p$-adic meromorphic function.*

In [Dw2-5], Dwork showed that Conjecture 6.3 is true for several examples. This includes the universal family of elliptic curves and a certain family of K3-surfaces. His idea is to reduce the problem to the classical overconvergent setting by establishing the existence of the so-called excellent lifting. The excellent lifting however rarely exists and thus one cannot hope that it would work in general. Another approach studied in [DS] and [W3] is to try to relax the overconvergent condition by using the weaker $c\log$-convergent condition. This approach pushes Dwork's trace formula to its full potential but again it cannot work in general by the counter-example in [W3]. The general case of Dwork's conjecture is being proved in our recent series of papers [W8-10] by introducing an entirely new method building on the fundamental Dwork-Monsky trace formula [Mo]. At this point, the abstract version of Dwork's conjecture for F-crystals and $\sigma$-modules, as proved in [W8-10], implies Conjecture 6.3 whenever the finiteness of the relative rigid cohomology with compact support is known such as the smooth proper liftable case. Our strategy to handle singular and open family is to directly work with the infinite rank version. This requires additional work which is in progress.

Thus, the $k$-th power slope $s$ L-function $L^{[k]}(s,f,T)$ is always a good function. An immediate application is the existence of a formula for the $k$-th moment slope $s$ sequence $S_{k,d}(s,f)$. For each $k$ and each $s$, we can write
$$L^{[k]}(s,f,T) = \frac{\prod_i (1-\gamma_i(k,s)T)}{\prod_j (1-\delta_j(k,s)T)},$$

where the product is finite in the complex absolute value case and infinite in the $p$-adic case. We have

$$\lim_{i\to\infty} \gamma_i(k,s) = 0 = \lim_{j\to\infty} \delta_j(k,s).$$

The formula for the $k$-th moment slope $s$ sequence is

$$S_{k,d}(s,f) = \delta_1(k,s)^d + \delta_2(k,s)^d + \cdots - \gamma_1(k,s)^d - \gamma_2(k,s)^d - \cdots, \ d = 1, 2, \cdots.$$

Such formulas imply the existence of a suitable equi-distribution theorem for the zeros and poles of the family of pure rational functions $Z_s(Y_x, T)$ parametrized by $x \in X$. In the complex absolute case, this is related to but different from the approach of Deligne-Katz via monodromy representations. In the $p$-adic absolute value case, such a result is completely new.

Once we know that $L^{[k]}(s, f, T)$ is a good function, we can then ask for the various RHs for $L^{[k]}(s, f, T)$. In the case of Theorem 6.2, there are three types of RHs (the complex, the $\ell$-adic and the $p$-adic). The situation is similar to Section 5. The L-function in Theorem 6.2 has integer coefficients in view of Theorem 3.8.

In the case of Conjecture 6.3, one can only ask the $p$-adic RH since the zeros and poles are $p$-adic numbers. This $p$-adic RH for the L-function in Conjecture 6.3 is apparently more difficult than the already mysterious $p$-adic RH for the L-function in Theorem 6.2, because there are now infinitely many zeros and poles in the case of Conjecture 6.3. In fact, we have no ideas what the general conjectural formulation for the $p$-adic RH should be. It seems un-reasonable to expect a finite upper bound for the denominators of the slopes of the zeros and poles of $L^{[k]}(s, f, T)$, although such a finite bound is conjectured to be true [W6] in the simplest elliptic family case.

A harder higher dimensional example is given by the following family of $n$-dimensional Calabi-Yau varieties:

$$X_1^{n+2} + X_2^{n+2} + \cdots + X_{n+2}^{n+2} + \lambda X_1 X_2 \cdots X_{n+2} = 0.$$

In the case $n = 1$ (elliptic curves), the $k$-th power slope $s$ ($s = 0$) L-function $L^{[k]}(s, f, T)$ in the $p$-adic case already contains a great deal of arithmetic information about classical and $p$-adic modular forms. In the case $n = 2$ (K3-surfaces), Conjecture 6.3 was proved by Dwork but its arithmetic consequences have not been explored. In the higher dimensional case $n > 2$, Conjecture 6.3 follows from our recent work. It would be of great interest to explore its arithmetic consequences. For instance, the $k$-th power slope $s$ ($s = 0$) L-function $L^{[k]}(s, f, T)$ in the $p$-adic case should be closely related to the arithmetic of the mysterious mirror map in mirror symmetry, see Lian and Yau [LY] for another relation bewteen Dwork's work and the mirror map. In some cases, the special value of the L-function $L^{[k]}(s, f, T)$ at $T = 1$ seems to be related to the conjectural $p$-adic L-functions of algebraic varieties defined over a number field. Very little is known in this direction.

# 7 Zeta functions of algebraic cycles

In this final section, following [W1], we define the zeta functions of algebraic cycles of a projective variety embedded in a given projective space. Several standard conjectures associated with such zeta functions are described.

We first recall the definition of the degree of a projective variety. Let $n$ and $m$ be positive integers with $n \leq m$. Let $X$ be a closed $n$ dimensional subscheme of the $m$ dimensional projective space $\mathbf{P}^m$ over $\mathbf{F}_q$ (the embedding will be fixed). We shall always work over the ground field $\mathbf{F}_q$. Let

$$S(X) = \oplus_{d=0}^{\infty} S_d(X)$$

be the homogeneous coordinate ring of $X$, where $S_d(X)$ consists of the homogeneous elements of degree $d$ and $S(X)$ is finitely generated by $S_1(X)$ as an $\mathbf{F}_q$-algebra ($S_0(X) = \mathbf{F}_q$). Let

$$l(d) = \dim_{\mathbf{F}_q} S_d(X).$$

For large $d$, $l(d)$ is a polynomial in $d$ of degree $n$, called the Hilbert polynomial of $X$. The degree of $X$, denoted by $\deg X$, is defined to be the coefficient of $d^n/n!$ in this polynomial. If $X$ is irreducible, an alternative definition of $\deg X$ is the intersection number $X \cdot H^n$, where $H$ is a hyperplane in $\mathbf{P}^m$. This means that the degree of $X$ is the number of intersection points cut out by the intersection of $n$ sufficiently general hyperplanes.

Fix a closed $n$ dimensional subscheme $X$ of $\mathbf{P}^m$ over $\mathbf{F}_q$. A prime cycle on $X/\mathbf{F}_q$ is a closed integral subscheme of $X$ (i.e., a reduced and irreducible closed subscheme). For each integer $0 \leq r \leq n$, we define the zeta function of $r$-cycles on $X$ to be the following formal power series

$$Z_{r-\text{cycles}}(X, T) = \prod_P (1 - T^{\deg P})^{-1},$$

where $P$ runs over all prime cycles of dimension $r$ on $X$ and $\deg P$ is the degree of $P$ viewed as a closed subscheme of $\mathbf{P}^m$. Since the definition of the degree of a variety depends on the embedding, thus the zeta function of $r$-cycles for $0 < r < n$ depends on the embedding of $X$ in the projective space $\mathbf{P}^m$.

Let $N_r(d)$ (resp. $M_r(d)$) be the number of prime $r$-cycles (resp. effective $r$ cycles) of degree $d$ on $X$. By a theorem of Chow and van der Waerden, the set of effective $r$-cycles of degree $d$ is parametrized (one-to-one) by an algebraic set (the Chow variety) in a projective space. Thus, $M_r(d)$ and $N_r(d)$ are all finite. The zeta function of $r$-cycles on $X/\mathbf{F}_q$ is well defined. Since $Z_{r-\text{cycles}}(X, T)$ is a power series with integer coefficients, it is trivially $p$-adic analytic in the open unit disk.

Alternative expressions for $Z_{r-\text{cycles}}(X, T)$ are given by

$$Z_{r-\text{cycles}}(X,T) = \prod_{d=1}^{\infty}(1-T^d)^{-N_r(d)}$$
$$= \sum_{d=0}^{\infty} M_r(d) T^d$$
$$= \exp(\sum_{d=1}^{\infty} \frac{W_r(d)}{d} T^d),$$

where $W_r(d)$ is the following weighted sum (each prime $r$-cycle of degree $k$ is counted $k$ times)
$$W_r(d) = \sum_{k|d} k N_r(k).$$

Geometrically, a point counted in $N_r(d)$ consists of $d$ conjugate "non closed points" of dimension $r$. Thus, $W_r(d)$ can be interpreted as the number of "non closed points" of degree $d$ and dimension $r$. In particular, $W_0(d)$ is just the number of $\mathbf{F}_{q^d}$-rational points on $X$ and $Z_{0-\text{cycles}}(X,T)$ is the classical zeta function of a projective scheme over a finite field. The zeta function $Z_{0-\text{cycles}}(X,T)$ of zero cycles is a rational function by Dwork's theorem. It is easy to show that for $0 < r < n$, the zeta function $Z_{r-\text{cycles}}(X,T)$ of $r$-cycles is never rational and in fact never complex meromorphic, as the radius of convergence for $Z_{r-\text{cycles}}(X,T)$ is zero as a complex function. However, some simple examples [W1] suggest that $Z_{r-\text{cycles}}(X,T)$ might be $p$-adic meromorphic in some cases. If $Z_{r-\text{cycles}}(X,T)$ is indeed $p$-adic meromorphic, it would imply that there are two sequences of $p$-adic integers $\alpha_i$ and $\beta_i$ approaching zero such that for all positive integers $d$, we have
$$W_r(d) = \beta_1^d + \beta_2^d + \cdots - \alpha_1^d - \alpha_2^d - \cdots.$$

This would generalize the classical formula for the number $W_0(d)$, in which there are only finitely many terms, because the function $Z_{0-\text{cycles}}(X,T)$ of zero cycles is rational.

To formulate our conjectures, we let $A_r(X)$ be the Chow group of $r$-cycles on $X/\mathbf{F}_q$ modulo rational equivalence. It is conjectured that $A_r(X)$ is a finitely generated abelian group. Let $A_r^+(X)$ be the monoid in $A_r(X)$ generated by the effective $r$-cycles on $X/\mathbf{F}_q$. This is the effective cone in $A_r(X)$. The monoid $A_r^+(X)$ is in general not a finitely generated monoid for $0 < r < n$. Although we do not have a proven counter-example, we feel that the $p$-adic meromorphic continuation of $Z_{r-\text{cycles}}(X,t)$ might be false if the monoid $A_r^+(X)$ is not finitely generated. This is a little inconsistent with the general philosophy that L-functions arising in a natural way from algebraic geometry are analytically good functions. Perhaps, one should not push the general principle too far. Thus, all our conjectures assume that $A_r^+(X)$ is a finitely generated monoid. In fact, we are confident about the truth of our conjectures only under the stronger assumption that $A_r(X)$ is of rank one. However, we will state our conjectures under the weaker assumption that the monoid $A_r^+(X)$ is finitely generated. The meromorphic continuation conjecture is

**Conjecture 7.1 (meromorphic continuation).** *Let $X$ be an $n$ dimensional projective variety in $\mathbf{P}^m$ over $\mathbf{F}_q$. Assume that the monoid $A_r^+(X)$ is a finitely generated monoid. Then $Z_{r-\text{cycles}}(X,T)$ is a $p$-adic meromorphic function.*

Since the zeta function of $r$-cycles is supposed to be a $p$-adic meromorphic function, we can only ask for its $p$-adic RH. Any RH is really a finiteness property about zeros and poles. Thus, the $p$-adic RH in our current situation should be a finitness property on the slopes of the zeros and poles. Since there are infinitely many zeros and poles in general for the zeta function of $r$-cycles, the set of slopes for the zeros and poles is unbounded. Thus, the best we can hope for is to bound the denominators of the slopes. In this direction, we propose

**Conjecture 7.2 ($p$-adic RH).** *Let $X$ be an $n$ dimensional projective variety in $\mathbf{P}^m$ over $\mathbf{F}_q$. Assume that the monoid $A_r^+(X)$ is a finitely generated monoid. Then the denominators of the slopes (which are rational numbers) of the reciprocal zeros and reciprocal poles of $Z_{r-\text{cycles}}(X,T)$ are bounded by a constant which may depend on $X$.*

In connection with Tate's conjecture about the order of poles, we propose

**Conjecture 7.3 (order of pole).** *Let $X$ be an $n$-dimensional smooth projective variety in $\mathbf{P}^m$ over $\mathbf{F}_q$. Assume that the monoid $A_r^+(X)$ is a finitely generated monoid. Then, the rank of $A_r(X)$ is equal to the order of pole of $Z_{r-\text{cycles}}(X,T)$ at $T=1$.*

All three conjectures are true in the two extremal cases $r=0,n$. The first accessible new case is the case of divisors, i.e., when $r=n-1$. If $r=n-1$, the above three conjectures are known to be true in the case when $A_{n-1}(X)$ has rank one [W1]. In this case, Conjecture 7.2 was not stated and hence not proved in [W1], but its truth follows from the meromorphy proof given there. Furthermore, Conjecture 7.3 is always true in the divisor case $r=n-1$ without the rank one assumption.

Thus, even in the case of divisors, Conjectures 7.1-7.2 are not known to be true in general if $A_{n-1}(X)$ has rank greater than 1. At this point, the most fundamental conjecture seems to be Conjecture 7.1. Once Conjecture 7.1 is proved (if it is true), its proof should give a great deal of information about the other two conjectures. But no single example is known for $1 \leq r \leq n-2$. The simplest substantial example is to consider $Z_{1-\text{cycles}}(\mathbf{P}^3,T)$ (counting space curves of degree $d$ in $\mathbf{P}^3$ as $d$ varies), which seems already sufficiently difficult and requires new mathematics.

Tate's conjecture says that the order of pole at $T = q^{-r}$ of the zeta function of zero cycles is equal to the rank of the group of algebraic $r$-cycles modulo $\ell$-adic homological equivalence. For a variety over a finite field, this latter rank should equal to be the rank of the Chow group $A_r(X)$. With this modification, Tate's conjecture can be reformulated as

**Conjecture 7.4 (Tate).** *Let $X$ be an $n$-dimensional smooth projective variety over $\mathbf{F}_q$. Then, the rank of $A_r(X)$ is equal to the order of pole of the zeta function $Z_{0-\text{cycles}}(X, T)$ at $T = q^{-r}$.*

As indicated in the above, Conjecture 7.3 is known to be true in the divisor case $r = n - 1$. This suggests that Tate's conjecture might also be provable in the divisor case $r = n - 1$ if the monoid $A_{n-1}^+(X)$ is finitely generated. Combining the previous two conjectures together, we obtain the following conjecture relating the zeta function of $r$-cycles to the zeta function of zero cycles.

**Conjecture 7.5.** *Let $X$ be an $n$-dimensional smooth projective variety in $\mathbf{P}^m$ over $\mathbf{F}_q$. Assume that the monoid $A_r^+(X)$ is a finitely generated monoid. Then, the order of pole of $Z_{r-\text{cycles}}(X, T)$ at $T = 1$ is equal to the order of pole of $Z_{0-\text{cycles}}(X, T)$ at $T = q^{-r}$.*

Based on the Griffith-Katz counter-example about 1-cycles on certain 3-fold, which applies only to varieties over fields with transcendental elements, we remarked in [W1] that one should not expect the equality between the order of pole at $T = 1$ of $Z_{r-\text{cycles}}(X, T)$ and the order of pole at $T = 1$ of $Z_{(n-r)-\text{cycles}}(X, T)$. This remark is misleading. As Sòule pointed out to me sometime ago, one should expect the equality of those two numbers by Beilinson's conjectures.

The special value at $T = 1$ of the zeta function of $r$-cycles is undoubtedly related to the torsion order of $A_r(X)$ and certain regulator of $A_r(X)$. This is proved in [W1] in the divisor case $r = n - 1$. We leave it to interested readers to find a conjectural formula for the special value for other $r$.